\documentclass[11pt]{article}

\usepackage[T1]{fontenc}
\usepackage[english]{babel}
\usepackage{amsmath}
\usepackage{amsfonts}
\usepackage{amssymb}
\usepackage{graphicx}
\usepackage{float}
\usepackage{multicol}
\usepackage{color}
\definecolor{gris25}{gray}{0.55}
\usepackage{shadow}
\usepackage{makeidx}
\usepackage{subfigure}
\begin{document}

\newcommand{\be}{\begin{equation}}
\newcommand{\ee}{\end{equation}}
\newcommand{\bd}{\begin{displaymath}}
\newcommand{\ed}{\end{displaymath}}
\newcommand{\lijn}{\vspace{3.5mm}}
\newcommand{\ba}{\begin{eqnarray}}
\newcommand{\ea}{\end{eqnarray}}
\newcommand{\ban}{\begin{eqnarray*}}
\newcommand{\ean}{\end{eqnarray*}}
\newcommand{\ul}{\underline}
\newcommand{\ol}{\overline}
\newcommand{\ra}{\rightarrow}
\newcommand{\LL} {I\!\!L}
\newcommand{\C} {I\!\!\!\!C}
\newcommand{\R} {\mathbb{R}}
\newcommand{\PP} {\mathbb{P}}
\newcommand{\E} {\mathbb{E}}
\newcommand{\N} {\mathbb{N}}
\newcommand{\I} {1\!\! 1}
\newcommand{\ZZ} {\mathbb{Z}}
\newcommand{\cov} {\mbox{cov}}
\newcommand{\var} {\mbox{Var}}
\newcommand{\IBS} {\mbox{IBS}}
\newcommand{\sign} {\mbox{sign}}
\newcommand{\ds}{\displaystyle}
\def\convlaw{\renewcommand{\arraystretch}{0.5}
\begin{array}[t]{c}\stackrel{{\cal D}}{\rightarrow} \\
{}\end{array}\renewcommand{\arraystretch}{1}}
\newcommand{\MotsCles}[1]{\par\noindent
{\textbf{Mots cl\'es}: #1}}

\newcommand{\Resume}[1]{\par\noindent
{\textbf{R\'esum\'e}\/: #1}}

\newcommand{\Keywords}[1]{\par\noindent
{\textbf{Keywords}\/: #1}}

\newcommand{\Abstract}[1]{\par\noindent
{\textbf{Abstract}\/: #1}}
\begin{center}
{\Large
	{\sc Fast change point analysis on the Hurst index of piecewise fractional Brownian motion
		 %Off-line detection of multiple change points in the Hurst parameter of piecewise fractional Brownian motion by the Filtered Derivative with p-Value method
	}
}
\bigskip

 Pierre, R.~BERTRAND${}^{1,2}$ 
 { \it Pierre.Bertrand@math.univ-bpclermont.fr}\\
 \medskip
 Mehdi FHIMA${}^{2}$ 
  { \it Mehdi.Fhima@math.univ-bpclermont.fr}\\
  \medskip
 Arnaud GUILLIN${}^{2}$
  { \it Arnaud.Guillin@math.univ-bpclermont.fr}\\
  
{\it
 ${}^{1}$   INRIA Saclay \\
${}^{2}$  Laboratoire de Math\'ematiques, UMR CNRS 6620\\
\& Universit\'e de Clermont-Ferrand II, France

}
\end{center}

\Resume{\small Dans cette pr\'esentation, nous introduisons une nouvelle m\'ethode de d\'etection de ruptures sur  l'indice de Hurst, pour un mouvement brownien fractionnaire par morceaux. En premier lieu, nous d\'efinissons le mod\`ele et le probl\`eme statistique. La m\'ethode propos\'ee est une transposition de la m\' ethode FDpV \`a l'estimation de l'indice de Hurst.
La m\' ethode FDpV (d\'eriv\'ee filtr\'ee avec \textit{p-value}) a \'et\'e  introduite pour d\'etecter des ruptures sur la  moyenne par Bertrand {\it et al.} (2011).  La statistique sous-jacente de la technologie FDpV est un nouvel estimateur  de l'indice de Hurst, appel\'e  statistique de Bernouilli des accroissements (IBS). \`A la fois les m\'ethodes FDpV et IBS ont une complexit\'e lin\'eaire par rapport \`a la taille de la s\'erie d'observation, aussi bien en temps de calcul que pour la m\'emoire, donc \'egalement leur combinaison.}
\MotsCles{D\'etection de ruptures, D\'eriv\'ee Filtr\'ee, mouvement Brownien fractionnaire par morceaux, param\`etre de Hurst, Satistique de Bernouili des accroissements.}\\
\noindent
\Abstract{\small 
In this presentation, we introduce a new method for change point analysis on the Hurst index for a piecewise fractional Brownian motion. We first set the model and the statistical problem. The proposed method is a transposition of the FDpV (Filtered Derivative with p-value) method introduced for the detection of change points on the mean in Bertrand \textit{et al.} (2011) to the case of changes on the  Hurst index. The underlying statistics  of the FDpV technology is a new statistic estimator for Hurst index, so-called  Increment Bernoulli Statistic (IBS). Both FDpV and IBS are methods with linear time and memory complexity, with respect to the size of the series. Thus the resulting method for change point analysis on Hurst index reaches also a linear complexity. 
}

\section*{Introduction}
Recent measurement methods allow us to record and to stock large data sets, so called "the data deluge".   
For instance,  today technology allows recording of heartbeat series during  24 hours leading to data sets of size $n\ge
100,000$, and very high frequency (VHF) financial series leads to data size $n \geq 40,000$. 
Tomorrow, many other series will be recorded at VHF leading  to millions of data. 

Large or huge series with small meshes of time can be described as continuous time processes observed at discrete times.
Such a stochastic process $X$ belongs to a certain class of model, that is  $\ds X\in \mathcal{M}=\{X_{\theta},\; \theta\in \Theta\}$, where $ \Theta$ is a subset of $\R^d$ and $d$ is the dimension of the model.
The structural parameter $\theta$ is believed to provide relevant information on the system which generates the series, and statisticians have to estimate it. 

A slightly different approach is based on change point analysis: 
The structural parameter $\theta$ is assumed to be piecewise constant with an unknown configuration of change ${\bf \tau}$.  In this framework, the first task of statisticians is the estimation of the location of the change points and the second one could be the estimation of the structural parameters between change points. There is a huge literature on change point analysis and model selection  since the fifties, see  e.g. Basseville and  Nikiforov~(1993),
Brodsky and Darkhovsky~(1993), Cs\"{o}rgo and Horv\'ath~(1997), Birg\'e and Massart~(2007)  or Bertrand {\it et al.}~(2011) and the references therein. However, most of the studies are devoted to change on the mean, on the variance or on the regression parameters. 
But relevant informations are also provided by the time dependence structure of the process, see e.g. Ayache and Bertrand~(2011) and Khalfa {\it et al.}~(2011).  
Fractional Brownian motion (fBm) is a paradigmatic example of such process, indeed fBm  is a zero mean Gaussian process depending on two parameters:  The Hurst index $H$  linked to the time structure and a scale parameter $\sigma$. 

In this presentation, we consider a simple model, that is a process $X$ which is a piecewise fBm with an unknown configuration of changes. Moreover, we set us in the frame of huge datasets,  and we focus our attention on time and memory complexity. These two reasons have lead us to propose a new change point procedure for detection of change on the Hurst index for a piecewise fBm. Our new procedure is the combination, on the one hand of the FDpV  method, introduced in Bertrand {\it et al.}~(2011) for fast and light detection of change on the mean, variance or regression parameter, and on the other hand of the Increment Bernoulli Statistic (IBS) a new estimator for Hurst index, which is a variation on the Increment Ratio Statistic (IRS) estimator introduced in Bardet and Surgailis~(2010).

The rest of this paper is organized as follows:  At first, in Section~\ref{sec1}, we define  our  model of piecewise fBm. 
Next, in Section~\ref{sec2}, we introduce a new fast and robust estimator of the Hurst index of fBm, namely the Increment Bernoulli Statistic. Then, in Section~\ref{sec3}, we describe the transposition of the  FDp-V method to Hurst index. 

\section{Our model}
\label{sec1}
We observe a process $X$ at the discrete and regularly spaced time $t_i=i/n$, where $i=0,\dots,n$.
We assume the existence of a segmentation  $\tau=(\tau_k)_{k=0,\dots,K+1}$,  with $0=\tau_0<\tau_1<\dots<\tau_K<\tau_{K+1}=n$, such that the restriction of the process $X$ on each interval $(\tau_k,\tau_{k+1})$ for $k=0,\dots,K$ is a fBm with Hurst index $H_k$ and scale parameter $\sigma_k$. 
The integer $K$
corresponds to the number of change points and $(K+1)$ to the
number of segments. Stress that $K$ can be zero and in this case the process $X$ is a fBm.

Let us precise that in our roadmap the process $X$ should be almost surely with continuous paths. For this reason, the so-called piecewise fBm can not be defined by plugging a piecewise Hurst index into one of the representations of the fBm. Indeed, by doing so, the process $X$ would  almost surely have discontinuity at at each change point on Hurst index, and method for detecting change on the mean will also be efficient for detecting change on Hurst index.  Let us refer to  Taqqu and Samorodinitsky~(1994) as a reference book on the different representations of fBm and to Ayache and Taqqu~(2005) for the construction of multi-fBm by plugging a continuous time varying Hurst index into one of the fBm representations. A rather complicated solution to avoid the drawback of pathwise discontinuities due to Hurst index discontinuities, was proposed  in Benassi {\it et al.}~(2000) and cosigned by the first author. However, as point out by Antoine Ayache during private conversations held in 2004, for statistical applications, it suffices to cancel  the discontinuity by adding a correction term. The same solution is also adopted in  Bardet and Kammoun~(2008) .

The model having being specified, we are concerned with change point analysis on the Hurst parameter, where the number of change $K$ is unknown. There are few references on this problem. To our best knowledge, the only reference are Benassi {\it et al.}~(2000) and Bardet and Kammoun~(2008).

\section{Increment Bernoulli Statistic for fBm}
\label{sec2}
In this section, we investigate the properties of a new  estimator of the Hurst parameter of fBm, namely the Increment Bernoulli Statistic (IBS). IBS is a variation on IRS which has been introduced by Surgailis {\it et al.}~(2008) and applied to fBm by Bardet and Surgailis~(2010). Both IRS and IBS are fast and robust estimator of the Hurst index. By fast we mean estimator with linear time complexity, and by robust we mean estimator with invariant scaling property. The choice of the IBS instead of the IRS is motivated by the fact that the IBS is bit less expensive, in terms of time complexity, than the IRS. 

In the next section,  the IBS is used as the underlying estimator for the FDpV method, see $\eqref{def:D}$. For this reason, we define IBS for a every process $X$, even if we apply it to fBm in this section.  Let $X$ be a process observed at a family of   discrete times $t_k$, 
we define  the  second order increments by
$$\ds 
 \Delta X(t_k) := X(t_{k+2})-2X(t_{k+1})+ X(t_{k}).$$ 
 Then, the Increment Bernoulli Statistic (IBS) is based on the comparison of the signs of consecutive second order increments. The results of these comparisons will be equal to $1$ if the consecutive second order increments have the same sign, and $0$ otherwise. Hence, this explains the name of our new estimator, that is to say: Increment Bernoulli Statistic (IBS) which~is~given~by 
$$ \IBS_{n}(X)= \frac{1}{n-2} \sum_{k=0}^{n-3} \psi \left(
\Delta X(t_k) , \Delta X(t_{k+1}) \right) \label{IBS} $$
where $\psi(\cdot,\cdot)$ is described as follows $\psi(x,y)=0$ if $\sign(x)=\sign(y)$ and $0$ otherwise, where $\sign(z)$ denotes the sign of $z$. 
Let us remark that IBS is scale invariant: Indeed, since  $\psi(\sigma x,\sigma y)= \psi(x,y)$ for $\sigma>0$, then $\IBS_{n}(\sigma X)=\IBS_{n}(X)$.

\begin{figure}[htbp]
\begin{center}
\begin{tabular}{c}
\includegraphics[width=10cm,height=4cm]{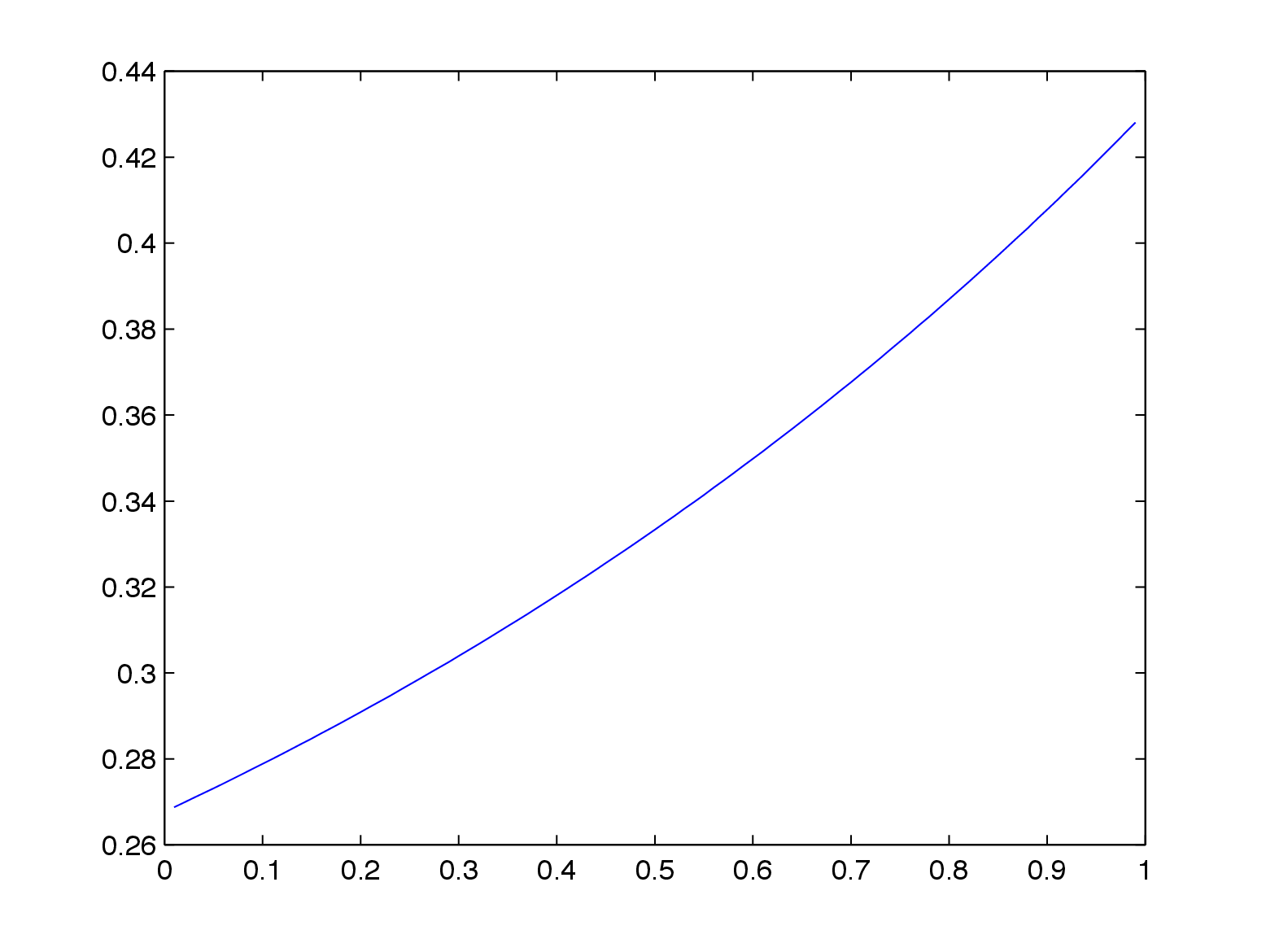}
\end{tabular}
\caption{\emph{The graph of $\Lambda(H)$.}} \label{fig1}
\end{center}
\end{figure}

When $X$ is a fBm, that is $X=B_H$ with a Hurst index $H \in (0,1)$, then 
the IBS converges in distribution to a continuous monotonic increasing function $\Lambda(H)$ defined as follows 
 \ban 
 \Lambda(H) &:=&\Lambda_0 \left(\rho(H) \right)  \\
 \Lambda_0(r) &:=&  \frac{1}{\pi} \arccos (-r) \\
 \rho (H) &=& \left (-3^{2H}+2^{2H+2}-7 \right)\left (8-2^{2H+1}\right)^{-1}
 \ean
where $ \rho(H) \in (-1,1)$ represents the correlation between two successive second order increments. The graph of $\Lambda(H)$ is given in Figure~\ref{fig1}. Then, due to the fact that $\Lambda(\cdot)$ is a reversible function, it is easy to deduce an estimator of the Hurst parameter $H$ given by $\widehat{H}_n = \Lambda^{-1}(\IBS_{n}(B_H))$.  Furthermore, we note that the function $\phi(\cdot,\cdot)=\psi(\cdot,\cdot)-\Lambda (H) $ is  a Hermite function with rank equal to 2.  Then, by  applying the Breuer-Major theorem, see {\it e.g.}  Arcones~(1994)[Theorem 4, p.2256] or Nourdin \emph{et al.}~(2010)[Theorem 1, p.2], we can deduce the following Central Limit Theorem (CLT): 
$$
\sqrt{n} \left( \IBS_{n}(B_H)-\Lambda(H) \right) \convlaw \mathcal{N}(0,\sigma^2(H)) \label{TCLfBm},
$$
 where the sign $\convlaw$ means convergence in distribution and the asymptotic variance $\sigma^2(H)$ is given by  
\ban \sigma^2(H)=\sum_{j \in \ZZ} \cov \left( \psi\left(
\Delta B_H(t_0) , \Delta B_H(t_{1}) \right), \psi \left(
\Delta B_H(t_j) , \Delta B_H(t_{j+1}) \right) \right). \ean 
The main advantages of the IBS method are primarily its efficiency in terms of time and memory complexity, and secondarily its robustness with respect to  scaling properties of the fBm. At first, we calculate by recurrence the second order increments $\left ( \Delta B_H(t_k) \right )_{0 \leq k \leq n-3}$. This first step is performed in time and memory complexity on $\mathcal{O}\left(n\right)$. Next, the computing of $\IBS_{n}(B_H)$ requires roughly $n$ tests + $n \times p_a(H)$ additions + 1 division, where $p_{a}(H) =\Lambda(H) \in (0,1)$ corresponds to the probability that two consecutive second order increments have the same sign. Then, we apply the 
Newton algorithm to compute the inverse of the function $\Lambda(\cdot)$. Moreover, we note that the function $\psi(\cdot,\cdot)$ satisfy the scale invariant property, {\it i.e.} for all $C \in \R$, $\psi(C.X, C.Y)=\psi(X,Y)$. This means that the multiplication of $B_H$ by any scaling coefficient $C$ does not  impact the estimation of the Hurst index, since $\IBS_{n}(B_H)=\IBS_{n}(C.B_H)$. Hence, this proves the robustness of the IBS method with respect to scaling. %properties of the fBm.
\section{Filtered Derivative with p-Value method}
\label{sec3}
In this section, we describe the Filtered Derivative with p-value method (FDp-V). First, we define the Filtered Derivative function. Next, we describe the two steps of the FDp-V method: Step~1 is based on Filtered Derivative and select the potential change points, whereas Step~2 calculate the p-value associated to each potential
change point, for disentangling right change points and false alarms.

 At first, we note that $\Lambda(H)$ is a continuous monotonic increasing function of $H$, see Figure~\ref{fig1}. So, the detecting of change points  on the Hurst parameter $H$ is equivalent to detecting   change points on  $\Lambda(H)$.  Consequently, the estimator  $\IBS_{n}(B_{H})$ of the parameter $\Lambda(H)$ is used as the underlying estimator for the FDpV method. We reefer to Bertrand {\it et al.}~(2011) and Bertrand and Fhima~(2009) for the introduction of FDpV technology and its numerical efficiency.
 Let us stress that  the choice of the direct estimator $\ds\widehat{H}_n = \Lambda^{-1}(\IBS_{n}(B_H))$ as the underlying estimator for the FDpV would be more expensive in term of numerical complexity.
\subsubsection*{Filtered Derivative function}
Let X be a piecewise fBm observed at a family of discrete times $t_j=j/n$, for $j=0,\dots,n$. 
The Filtered Derivative for IBS is defined as the difference between the estimators of the parameter $\Lambda(H)$ computed on two sliding windows respectively at the right and at the left of the index $k$, both of size $A$, that is specified by the  following function 
\ba
D(k,A)  = \IBS \left(X,k,A\right)-\IBS \left(X,k -A,A\right) \text{ for } k \in [A,n-A]
\label{def:D}
\ea
where  $$ \IBS \left(X,k,A\right)= A^{-1} \sum_{j=k+1}^{k+A} \psi \left( \Delta X(t_k) , \Delta X(t_{k+1}) \right)$$ is an estimator of $\Lambda(H)$ on the sliding box $[k+1,\,k+A]$. It is easy to see that the Filtered Derivative function $D$ is computed by recurrence with linear time and memory complexity. Eventually, this method consists on filtering data by computing the estimators of the parameter $\Lambda(H)$ before  applying a discrete derivation. This construction explains the name  given by Benveniste and Basseville~(1984), the so-called Filtered Derivative method. 
\subsubsection*{Step~1: Detection of potential change points}
\begin{figure}[htbp]
\begin{center}
\begin{tabular}{c}
\includegraphics[width=12cm,height=5.4cm]{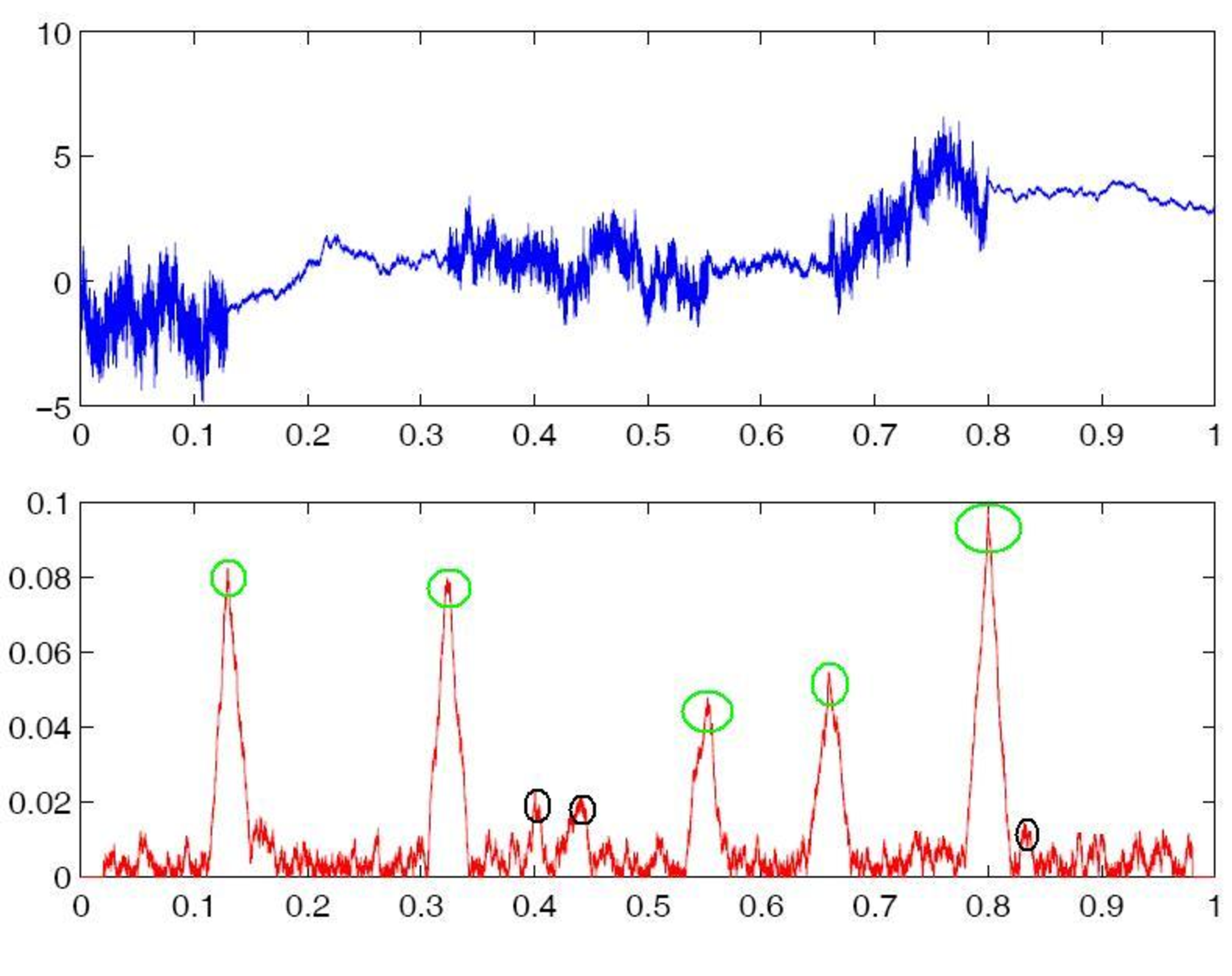}
\end{tabular}
\caption{\emph{Detection of potential change points. Above: Simulated piecewise fBm with five change points in the Hurst parameter. Below: Filtered Derivative function |D|.}} \label{fig2}
\end{center}
\end{figure}

In order to detect the potential change points, we test the null hypothesis $(\mathcal{H}_0)$ of no change in the Hurst parameter $H$
% $$(\mathcal{H}_0): \    H_1=H_2=\cdots=H_{n-1}=H_n$$
against the alternative hypothesis $(\mathcal{H}_1)$ indicating the existence of  at least one change point
$$(\mathcal{H}_1): \    \text{There is an integer } K \in \N^{*} \text{ and }
0=\tau_0<\tau_1<\dots<\tau_K<\tau_{K+1}=n \text{ such that }$$
$$
H_1 = \cdots = H_{\tau_{1}} \neq H_{\tau_{1}+1} =
\cdots = H_{\tau_{2}} \cdots \neq H_{\tau_{K}+1} =
\cdots = H_{\tau_{K+1}}.
$$
where $H_j=H_{\tau_{k}}$ is the value of the Hurst parameter at $t_j \in  [\tau_{k-1}/n,\tau_{k}/n)$.  \\
 Now, we fix a probability of type I error at level $p_1^*$,
and we determine the corresponding critical value $C_1$ given  by $$\PP \left( \max_{k \in [A:n-A]} |D(k,A)| >C_1| \mathcal{H}_0 \text{ is true} \right)=p_1^*.$$ 
Of course, such a probability is usually not available, so that we only consider the asymptotic distribution of the maximum of $|D|$. Then, the change points $\widetilde{\tau}_k$ is selected as a potential change point if its local maxima  satisfy $\ds |D(\widetilde{\tau}_k,A)|>C_1$. We remark through the graph of the function $|D|$ that there are not only the "right hats" (surrounded in green in Figure~\ref{fig2}) which gives the right change points, but also  false alarms (surrounded in black in Figure~\ref{fig2}). Consequently, we have introduced another step in order to keep just the right change points.
\subsubsection*{Step~2: Elimination of false alarms}
The list of potential change points $\displaystyle \left(\widetilde{\tau}_1,\ldots,\widetilde{\tau}_{K_{\max}} \right)$ obtained at step 1 contains right change points but also false alarms. In the second step, a test is carried out to remove the false alarms from the list of change points found at step 1. More precisely, for all potential change point $\widetilde{\tau}_{k}$, we test whether the Hurst parameter is the same on the two successive intervals $(\widetilde{\tau}_{k-1}/n,\widetilde{\tau}_{k}/n)$ and $(\widetilde{\tau}_{k}/n,\widetilde{\tau}_{k+1}/n)$, or not. Formally, for all $1 \le k \le K_{\max}$, we apply  the following hypothesis testing  $$  (\mathcal{H}_{0,k}):   H_k =   H_{k+1} \quad\text{versus}\quad  (\mathcal{H}_{1,k}): H_k \neq   H_{k+1},$$ 
where $H_k$ is the value of $H$ on the
segment $(\widetilde{\tau}_{k-1}/n,\widetilde{\tau}_{k}/n)$. By
using this second test, we calculate new p-values $\ds
(\widetilde{p}_{1},\ldots,\widetilde{p}_{K_{\max}})$ associated
respectively to each potential change points $\ds
(\widetilde{\tau}_1,\ldots,\widetilde{\tau}_{K_{\max}})$. Then, we
only keep the change points which have a p-value smaller than
a critical level denoted $p_2^{*}$. By doing so, we obtain a subset $\displaystyle \left(\widehat{\tau}_{1},\ldots,\widehat{\tau}_{\widehat{K}} \right)$ of the first list which represents the estimators of the change points in the Hurst parameter of mBm.

\section*{Conclusion}
In conclusion, it appears that the combination of the FDpV and the IBS methods provides a fast (time) and cheap (memory) algorithm to the detection of change points on the Hurst parameter of  piecewise fBm. So, this algorithm is adapted to segment random signals with large datasets. In future work, we will develop the FDpV $+$ IBS method in order to detect abrupt changes on parameters of real data drawn from financial and physiological domains.\\ 
\newpage
\noindent {\large{\bf References}}
\par
\medskip
\noindent [1] Antoch, J. and Hu{\v{s}}kov{\'a}, M. (1994)  Procedures for the detection of multiple changes in series of independent observations, {\it Asymptotic statistics ({P}rague, 1993)}, Contrib. Statist., 3--20. Physica, Heidelberg.

\noindent [2] Arcones, M. A (1994) Limit theorems for nonlinear functionals of a stationary Gaussian sequence of vectors, {\it Ann. Probab}, 22, 2242--2274.

\noindent [3] Ayache, A. and Bertrand, P.R. (2011) Discretization error of wavelet coefficient for fractal like process, to appear in {\it Advances in Pure and Applied Mathematics}.

\noindent [4] Ayache, A. and Taqqu, M. S. (2005) Multifractional process with random exponent, {\it Publ. Mat}, 49, 459--486.

\noindent [5] Bardet, J. M. and Kammoun, I. (2008). Detecting abrupt changes of the long-range dependence or the self-similarity of a Gaussian process, {\it C. R. Math. Acad. Sci}. Paris, 346, 889--894.

\noindent [6] Bardet, J. M. and Surgailis, D. (2010) Measuring Roughness of Random Paths by Increment Ratios, to appear in {\it Bernoulli}.

\noindent [7] Basseville, M. and  Nikiforov, I. V. (1993) {\it Detection of abrupt changes: theory and application}, Prentice Hall Inc, Englewood Cliffs, NJ.  

\noindent [8] Benassi, A., Bertrand, P. R., Cohen, S., Istas, J. (2000) Identification of the Hurst exponent of a Step Multifractional Brownian motion, {\it Statistical Inference for Stochastic Processes}, 3, 101--111.

\noindent [9] Benveniste, A. and Basseville, M. (1984) Detection of abrupt changes in signals and dynamical systems: some statistical aspects, {\it Analysis and optimization of systems, {P}art~1 ({N}ice, 1984), volume 62 of Lecture Notes in Control and Inform. Sci.}, 145--155, Springer, Berlin.

\noindent [10] Bertrand, P. R. and M. Fhima, M. (2009) Filtered Derivative with p-Value Method for Multiple Change-Points Detection, {\it Proceeding of the 2nd International Workshop in Sequential Methodologies}.

\noindent [11] Bertrand, P. R., Fhima, M. and Guillin, A. (2011) Off-line detection of multiple change points by the Filtered Derivative with p-Value method, to appear in {\it Sequential Analysis}.

\noindent [12] Birg{\'e}, L. and Massart, P. (2007) Minimal penalties for {G}aussian model selection, {\it Probab. Theory Related Fields}, 138(1-2), 33--73.

\noindent [13] Brodsky, B. E. and Darkhovsky, B. S. (1993) {\it Nonparametric methods in change-point problems, volume 243 of Mathematics and its applications}, Kluwer Academic Publishers Group, Dordrecht. 

\noindent [14] Cs\"{o}rgo, M. and Horv\'ath, L. (1997) {\it Limit Theorem in Change-Point Analysis}, J. Wiley, New York.

\noindent [15] Khalfa, N., Bertrand, P. R., Boudet, G., Chamoux, A. and Billat, V. (2011), Heart Rate Regulation processed through wavelet analysis and change detection. Some case studies, {\it Submitted}.

\noindent [16] Nourdin, I., Peccati, G. and Podolskij, M. (2010) Quantitative Breuer-Major Theorems, HAL : hal-00484096, version 2.

\noindent [17] Samorodnitsky, G. and Taqqu, M. S. (1994) {\it Stable non-Gaussian random processes},  Chapman \& Hall.

\noindent [18] Surgailis, D., Teyssière, G. and Vai\v{c}iulis, M. (2008) The increment ratio statistic, {\it J. Multivariate Anal}, 99, 510--541. 

\end{document}